\documentclass[12pt]{article}
\usepackage{amsfonts}
\usepackage{amsmath}
\usepackage{enumitem}
\usepackage{sectsty}
\usepackage{cite}
\usepackage[margin=1in]{geometry}
\usepackage{amsthm}
\usepackage{mathtools}
\usepackage{amssymb}
\usepackage{graphicx}

\makeatletter
\def\thm@space@setup{%
  \thm@preskip=0.5cm plus 0.5cm 
  \thm@postskip=\thm@preskip 
}
\makeatother

\sectionfont{\large}
\subsectionfont{\normalsize}
\subsubsectionfont{\small}

\newtheorem{thm}{Theorem}[section]

\newtheorem{prp}[thm]{Proposition}

\theoremstyle{definition}
\newtheorem{dfn}[thm]{Definition}

\theoremstyle{definition}

\theoremstyle{definition}
\newtheorem{exa}[thm]{Example}

\theoremstyle{definition}
\newtheorem{obs}[thm]{Remark}

\usepackage{titling}

\usepackage{lipsum}

\begin{document}

\title{Generalising the Wallis Product}
\author{Joshua W. E. Farrell}
\date{\today}
\maketitle

\begin{abstract}
In 1655, John Wallis whilst at the University of Oxford discovered the famous and beautiful formula for $\pi$, now known as Wallis' Product.\cite{ARTICLE:4}
$$\frac{\pi}{2} = \frac{2}{1}\cdot\frac{2}{3}\cdot\frac{4}{3}\cdot\frac{4}{5}\cdot\frac{6}{5}\cdot\frac{6}{7}\cdot\frac{8}{7}\cdot\frac{8}{9}\cdots$$

Since then, several analogous formulae have been discovered generalising the original. One more modern proof of the Wallis Product and its relatives directly uses the Gamma Function. This short paper will use similar techniques to understand certain related classes of infinite products. Almost all results within this paper are new findings made by myself; when I should be revising or completing assignment work I find myself always going back to this. 
\end{abstract}

\section{Preliminaries}

\begin{dfn} \label{D: Gamma Function}
The Gamma Function is a well known analytic continuation of the factorial function. It is the only function that satisfies $\Gamma(z+1)=z\Gamma(z)$, $\Gamma(1)=1$ and $\ln(\Gamma(z))$ convex (see Bohr-Mollerup Theorem \cite{article:9}). It is an entire function (everywhere analytic), hence by Weierstrass Factorisation Theorem we have,\begin{align*}
\Gamma(z) = \frac{e^{-\gamma z}}{z}\prod^{\infty}_{n=1} \left(1+\frac{z}{n} \right)^{-1} e^{z/n}
\end{align*} which is valid for all $z\in\mathbb{C}\setminus \mathbb{Z}_{\leq 0}$ \cite{article:8}. Some \textit{popular} formulae are the functional equation and Euler's reflection formula respectively \begin{align*}
 \Gamma(z+1)=z\Gamma(z) \quad\quad\quad\quad \Gamma(z,)\Gamma(1-z)=\frac{\pi}{\sin(\pi z)},\,\,\, \forall z\notin\mathbb{Z}.
\end{align*}
\end{dfn}

\begin{dfn} \label{D: Catalan's Constant}
Catalan's constant is denoted by $K$ rather than $\mathbf{G}$ to avoid confusion with the Barnes G-function and given by, \cite{article:12}
\begin{align*}
K := \sum_{k=0}^{\infty}\frac{(-1)^{k}}{(2k+1)^{2}}  \approx 0.915965594177...
\end{align*}  
\end{dfn}

\begin{dfn}\label{D: Barnes G-Function}
The Barnes G-Function is a generalisation of the Gamma Function and extension of the superfactorial, given by the functional equation $G(z+1)=\Gamma(z)G(z)$, $G(1)=1$ and $G^{(3)}(x)\geq 0$, similar to Bohr-Mollerup Theorem. Like the Gamma Function, the Weierstrass definition will be used throughout.
\begin{align*}
G(z+1) = (2\pi)^{z/2} \exp \left(-\frac{z+z^{2}(1+\gamma)}{2} \right)\prod^{\infty}_{n=1}\left\lbrace \left(1+\frac{z}{n}\right)^{n}\exp \left(\frac{z^{2}}{2n}-z \right)\right\rbrace
\end{align*}
Which is valid for all $z\in\mathbb{C}$. \cite{article:10}
\end{dfn}

\begin{dfn} \label{D: GK Constant}
The Glaisher-Kinkelin constant can be defined by a limit, or by it's relation to the Riemann Zeta Function \cite{article:12},
\begin{align*}
A := \lim_{n\to\infty}\frac{(2\pi)^{n/2}n^{n^{2}/2 - 1/12}e^{-3n^{2}/4 + 1/12}}{G(n+1)}= \exp\left(\frac{1}{12} - \zeta'(-1) \right) \approx 1.2824271291...\, .
\end{align*}  
\end{dfn}

\begin{dfn} \label{D: Multiple Gamma Function}
The Multiple Gamma Function is a further generalisation of the Gamma function, given by a functional equation and again a log-convexity property \cite{ARTICLE:2}, \begin{align*}
\Gamma_{n+1}(z+1) = \frac{\Gamma_{n+1}(z)}{\Gamma_{n}(z)}, \,\, \forall z\in\mathbb{C}, n\in\mathbb{N}
\end{align*} with additional criteria, $\Gamma_{n}(1)=1$ and $\Gamma_{1}(z)=\Gamma(z)$. It should also be noted that $G(z)=\Gamma_{2}(z)^{-1}$. Some literature use a generalised G-Function instead of Gamma Function which is defined very similarly as $G_{n+1}(z+1)=G_{n}(z)G_{n+1}(z)$, but the differences for this paper are unimportant as $\Gamma_{n}(z)=G_{n}(z)^{(-1)^{n+1}}$.
\end{dfn}

\section{The Gamma Function}

\begin{thm}
\label{T: Gamma Product}
Let $i\in\mathbb{N}_{+}$, $\lbrace a_{1}, ..., a_{i}, b_{1}, ..., b_{i} \rbrace \in\mathbb{C}\setminus \mathbb{Z}_{\leq 0}$ and $a_{1}+...+a_{i}=b_{1}+...+b_{i}$. Then,
\begin{align*}
\frac{\Gamma (b_{1})...\Gamma (b_{i})}{\Gamma (a_{1})...\Gamma (a_{i})} =\prod^{\infty}_{n=0} \frac{(n+a_{1})...(n+a_{i})}{(n+b_{1})...(n+b_{i})}.
\end{align*}
\emph{Products of this form will be known throughout as Type-I products.} 
\begin{proof}
See \cite{article:11}, the proof is very similar to the later more complicated proofs of Theorem \ref{T: II Product} and Theorem \ref{T: N product} hence I have not included it.
\end{proof}
\noindent  \emph{A clear corollary of this and the functional equation, Definition \ref{D: Gamma Function}, is the following formula which is valid when} $\lbrace a_{1}, ..., a_{i}, b_{1}, ..., b_{i} \rbrace \in\mathbb{C}\setminus \mathbb{Z}_{-}$, $a_{1}+...+a_{i}=b_{1}+...+b_{i}$,
\begin{align*}
\frac{\Gamma (b_{1}+1)...\Gamma (b_{i}+1)}{\Gamma (a_{1}+1)...\Gamma (a_{i}+1)} =\prod^{\infty}_{n=1} \frac{(n+a_{1})...(n+a_{i})}{(n+b_{1})...(n+b_{i})}.
\end{align*} \emph{These will be used interchangeably when referring to this theorem dependant on convenience/aesthetic appeal. The second will be specifically referred to as the alternate equation.}
\end{thm}

\begin{exa}
\label{E: I Wallis Product}
\begin{align*}
\prod_{n=1}^{\infty} \frac{(2n)^{2}}{(2n-1)(2n+1)} =  \frac{\pi}{2}
\end{align*}\begin{proof}
Using Theorem \ref{T: Gamma Product} (alternate equation) with $i=2$ and the following constants,
\begin{align*}
a_{1} &= a_{2} = 0, \quad\quad b_{1} = -\frac{1}{2}, \quad b_{2} = \frac{1}{2},
\end{align*} we have
\begin{align*}
\prod_{n=1}^{\infty} \frac{(2n)^{2}}{(2n-1)(2n+1)} = \prod_{n=1}^{\infty} \frac{n^{2}}{(n-1/2)(n+1/2)} =  \frac{\Gamma(1/2)\Gamma(3/2)}{\Gamma(1)^{2}}.
\end{align*} From Definition \ref{D: Gamma Function}: $\Gamma(1)=1$, the functional equation implies $\Gamma(3/2)=\Gamma(1/2)/2$ and the reflection formula implies $\Gamma(1/2)=\sqrt{\pi}$. Which when combined with the equation above gives the famous infinite product discovered by Wallis. \cite{ARTICLE:4}\cite{ARTICLE:5}
\end{proof}
\end{exa}

\begin{exa}
\label{E: I Golden Ratio}
\begin{align*}
\prod_{n=0}^{\infty} \frac{(30n+9)(30n+21)}{(30n+5)(30n+25)}  = \frac{1+\sqrt{5}}{2} = \varphi 
\end{align*} \begin{proof}
Using Theorem \ref{T: Gamma Product} with $i=2$ and the following constants,
\begin{align*}
a_{1} &= \frac{3}{10}, \quad a_{2} = \frac{7}{10}, \quad\quad b_{1} = \frac{1}{6}, \quad b_{2} = \frac{5}{6}.
\end{align*} Then manipulation continues as in Example \ref{E: I Wallis Product} using the formulae from Definition \ref{D: Gamma Function}.
\end{proof}
\end{exa}

\begin{obs}
See that if we choose $i=1$ then the product simplifies to $\prod_{n=1}^{\infty} 1 = 1$, a very uninteresting case. Also notice, choosing integer values for the $a_{k}$ and $b_{k}$ is fairly uninteresting as the products will always converge to some collection of factorials (resulting in some rational number). 
\end{obs}

\begin{prp} \label{P: Root Two}
$\forall k\in\mathbb{N}_{+},$
\begin{align*}
\prod^{\infty}_{n=0} \frac{(2^{k+1}n+2)(2^{k+1}n+2^{k+1}-2)}{(2^{k+1}n+1)(2^{k+1}n+2^{k+1}-1)} = \underbrace{\sqrt{2+\sqrt{2+\sqrt{2+...}}}}_{k} 
\end{align*}
\begin{proof}
Using Theorem \ref{T: Gamma Product} we can rewrite the product in terms of the Gamma Function.
\begin{align*}
\prod^{\infty}_{n=0} \frac{(2^{k+1}n+2)(2^{k+1}n+2^{k+1}-2)}{(2^{k+1}n+1)(2^{k+1}n+2^{k+1}-1)}&= \prod^{\infty}_{n=0} \frac{(n+2^{-k})(n+1-2^{-k})}{(n+2^{-(k+1)})(n+1-2^{-(k+1)})} \\
&= \frac{\Gamma(2^{-(k+1)})\Gamma(1-2^{-(k+1)})}{\Gamma(2^{-k})\Gamma(1-2^{-k})}.
\end{align*} Then using the reflection formula from Definition \ref{D: Gamma Function} and the double angle formula, this simplifies to
\begin{align*}
\frac{\Gamma(2^{-(k+1)})\Gamma(1-2^{-(k+1)})}{\Gamma(2^{-k})\Gamma(1-2^{-k})} &= \frac{\frac{\pi}{\sin(2^{-(k+1)}\pi)}}{\frac{\pi}{\sin(2^{-k)}\pi)}}
= \frac{\sin(2^{-k}\pi)}{\sin(2^{-(k+1)}\pi)} 
= 2\cos\left({\frac{\pi}{2^{k+1}}}\right).
\end{align*} Iterating the double angle formula for cosine,  \begin{align*}
2\cos\theta = \sqrt{2+2\cos 2\theta} = \sqrt{2+\sqrt{2+2\cos 4\theta }} = ...\, ;
\end{align*} from which, it easily follows that
\begin{align*}
 2\cos\left({\frac{\pi}{2^{k+1}}}\right)&= \underbrace{\sqrt{2+\sqrt{2+\sqrt{2+...}}}}_{k} .
\end{align*} 
\end{proof}
\end{prp}

\begin{exa}
Utilising Proposition \ref{P: Root Two}, for $k=1$, we retrieve the famous infinite product for $\sqrt{2}$ found by Catalan in 1875 \cite{article:7}. For $k=2$, we retrieve a formula found by Sondow and Yi in 2010 \cite{ARTICLE:3}. And for $k=3$,
 \begin{align*}
\prod^{\infty}_{n=0} \frac{(16n+2)(16n+14)}{(16n+1)(16n+15)}&= \sqrt{2+\sqrt{2+\sqrt{2}}}.
\end{align*}
\end{exa}

\begin{obs} \label{O: Gamma}
Due to the vast amount of formula available for computation of the Gamma Function, one can construct Type-I products to converge to many different collections of constants; we have already seen $\pi/2$, the golden ratio and an $n$ length nested radical of $2$'s. Two of the most useful formulae are Euler's reflection formula and $\Gamma\left(1/2\right)=\sqrt{\pi}$. Suddenly Example \ref{E: I Wallis Product} is much more obvious and when one knows $\cos (\pi /5) = \varphi /2$, Example \ref{E: I Golden Ratio} quickly follows. I will illustrate the ease of constructing products in the following Example.
\end{obs}

\begin{exa} \label{E: I General Rational Solutions}
\begin{align*}
\frac{p}{q}  = \prod_{n=1}^{\infty}\frac{n(qn+p)}{(n+1)(qn+p-q)}
\end{align*} \begin{proof}
Using Theorem \ref{T: Gamma Product} (alternate formula) with $i=2$ and the following constants,
\begin{align*}
a_{1} &= 0, \quad a_{2} = \frac{p}{q}, \quad\quad b_{1} = 1, \quad b_{2} = \frac{p-q}{q},
\end{align*} where $p,q\in\mathbb{N}$. Again, manipulation continues as in Example \ref{E: I Wallis Product} and Proposition \ref{P: Root Two}.
\end{proof}
\end{exa}

\begin{obs}
Due to Example \ref{E: I General Rational Solutions}, we see that not only can we express all positive rational numbers in the Type-I product form, but also we can express all rational multiples of $\pi$ in the Type-I form. As we can multiple convergent products, taking Example \ref{E: I Wallis Product} and multiplying it with Example \ref{E: I General Rational Solutions}, we have a product for $p\pi / 2q$ (which clearly can be any positive rational multiple of $\pi$). Furthermore, we can express any rational multiple of any power of $\pi$ as an infinite product, by the same logic.
\end{obs}

\section{The Barnes G-function}

\begin{thm}
\label{T: II Product}
Let $i\in\mathbb{N}_{+}$, $\lbrace a_{1}, ..., a_{i}, b_{1}, ..., b_{i} \rbrace \in\mathbb{C}\setminus\mathbb{Z}_{<0} $, $a_{1}+...+a_{i}=b_{1}+...+b_{i}$ and $a_{1}^{2}+...+a_{i}^{2}=b_{1}^{2}+...+b_{i}^{2}$. Then,
\begin{align*}
\prod^{\infty}_{n=1}\bigg[ \frac{(n+a_{1})...(n+a_{i})}{(n+b_{1})...(n+b_{i})}\bigg]^{n} = \frac{G(a_{1}+1)...G(a_{i}+1)}{G(b_{1}+1)...G(b_{i}+1)} .
\end{align*}
\emph{Products of this form will be known throughout as Type-II products.}
\begin{proof}
Let $\lbrace a_{1}, a_{2},...,a_{i}, b_{1}, b_{2},...,b_{i} \rbrace \in \mathbb{C}\setminus\mathbb{Z}_{<0}$ , then from the Weierstrass definition of the Barnes G-function (Definition \ref{D: Barnes G-Function}) we have, 
\begin{align*}
\frac{G(a_{1}+1)...G(a_{i}+1)}{G(b_{1}+1)...G(b_{i}+1)} &= (2\pi)^{T_{1}/2} e^{S}\prod^{\infty}_{n=1}\bigg\lbrace \frac{(1+a_{1}/n)^{n}...(1+a_{i}/n)^{n}}{(1+b_{1}/n)^{n}...(1+b_{i}/n)^{n}}\exp \Big(\frac{T_{2}}{2n}-T_{1} \Big)\bigg\rbrace,  \\ \\
\text{where}\quad S &= -\frac{1}{2}\big( T_{1}+T_{2}(1+\gamma) \big), \\
T_{1} &= (a_{1}+a_{2}+...+a_{i})-(b_{1}+b_{2}+...+b_{i}), \\
T_{2} &= (a_{1}^{2}+a_{2}^{2}+...+a_{i}^{2})-(b_{1}^{2}+b_{2}^{2}+...+b_{i}^{2}).
\end{align*}
If the constants satisfy the theorem criteria then, $T_{1}=0$ and $T_{2}=0$ (which further implies $S=0$). Simplifying the product gives the required formula. 
\end{proof}
\end{thm}

\begin{exa}\label{E: II Exponential Catalan}
\begin{align*}
\prod_{n=1}^{\infty}\bigg[\frac{(4n-1)^{3}(4n+3)}{(4n-3)(4n+1)^{3}}\bigg]^{n} &= e^{2K/\pi}
\end{align*} 
\begin{proof}
Using Theorem \ref{T: II Product} with $i=4$ and the following constants,
\begin{align*}
a_{1} &= a_{2} = a_{3} = -\frac{1}{4}, \quad a_{4} = \frac{3}{4}, \quad
b_{1} = -\frac{3}{4}, \quad b_{2} = b_{3} = b_{4} = \frac{1}{4}.
\end{align*}
Clearly, \begin{align*}
\prod_{n=1}^{\infty}\bigg[\frac{(4n-1)^{3}(4n+3)}{(4n-3)(4n+1)^{3}}\bigg]^{n} &= \prod_{n=1}^{\infty}\bigg[\frac{(n-1/4)^{3}(n+3/4)}{(n-3/4)(n+1/4)^{3}}\bigg]^{n} = \frac{G(3/4)^{3}G(7/4)}{G(1/4)G(5/4)^{3}} ,
\end{align*}
and it remains to prove that \begin{align*}
\frac{G(3/4)^{3}G(7/4)}{G(1/4)G(5/4)^{3}} = e^{2K/\pi}.
\end{align*}
Manipulation of the Barnes G-function is required, but this is not the purpose of the paper. I will conduct this proof, but others will be left out as the method is very similar. Firstly, we use the functional equation of the Barnes G-function to reduce all arguments such that they are between $0$ and $1$. \begin{align*}
\frac{G(3/4)^{3}G(7/4)}{G(1/4)G(5/4)^{3}} &= \frac{\Gamma(3/4)}{\Gamma(1/4)^{3}}\left(\frac{G(3/4)}{G(1/4)}\right)^{4}
\end{align*}
$G(3/4)$ and $G(1/4)$ can be written in terms of the Gamma Function as follows \cite{article:6},
\begin{align*}
G(1/4) &= A^{-9/8}\Gamma^{-3/4}(1/4)e^{3/32 - K/4\pi} \\
G(3/4) &= A^{-9/8}\Gamma^{-1/4}(3/4)e^{3/32 + K/4\pi} 
\end{align*} Hence,
\begin{align*}
\frac{G(3/4)}{G(1/4)} = \frac{\Gamma^{3/4}(1/4)}{\Gamma^{1/4}(3/4)}e^{K/2\pi}
\end{align*}
and the result easily follows.
\end{proof}
\end{exa}

\begin{obs}
As a side-note, I believe this to be a very important evaluation. The Dirichlet Beta Function, defined by \begin{align*}
\beta(s) = \sum_{k=0}^{\infty}\frac{(-1)^{k}}{(2k+1)^{s}}
\end{align*} for $\Re(s)>0$, has the specific evaluation $\beta'(-1) = 2K/\pi$  \cite{article:13}, that is the product in Example \ref{E: II Exponential Catalan} converges to $\exp(\beta'(-1))$. I have studied this in much greater detail and conjecture a general formula for all negative odd integers (there already exists a formula for positive integers). I have proved many lemmas relating to my conjecture and have strong evidence using computer tools such as PARI to support it further. I hope to have the full proof complete soon.
\end{obs}

\begin{exa}\label{E: II Rational}
\begin{align*}
\prod_{n=1}^{\infty} \left[\frac{(3n+1)(3n+5)^{2}(3n+7)}{(3n+2)(3n+4)^{2}(3n+8)}\right]^{n} &= \frac{4}{5}
\end{align*}
\begin{proof}
Using Theorem \ref{T: II Product} with $i=4$ and the following constants,
\begin{align*}
a_{1} = \frac{1}{3},\quad a_{2} = a_{3} = \frac{5}{3},\quad a_{4} = \frac{7}{3},\quad b_{1} = \frac{2}{3},\quad b_{2} = b_{3} = \frac{4}{3}, \quad b_{4} = \frac{8}{3}.
\end{align*}
The result follows from similar manipulation to Example \ref{E: II Exponential Catalan}.
\end{proof}
\end{exa}

\begin{prp}\label{P: Double Product}
Let $\lbrace a_{1}, ..., a_{i}, b_{1}, ..., b_{i} \rbrace \in\mathbb{C}$, $a_{1}+...+a_{i}=b_{1}+...+b_{i}$ and $a_{1}^{2}+...+a_{i}^{2}=b_{1}^{2}+...+b_{i}^{2}$. Then,
\begin{align*}
\prod^{\infty}_{n=1}\bigg[ \frac{(n+a_{1})...(n+a_{i})}{(n+b_{1})...(n+b_{i})}\bigg]^{n} = \prod_{r=0}^{\infty}\prod_{n=1}^{\infty}\frac{(n+r+a_{1})...(n+r+a_{i})}{(n+r+b_{1})...(n+r+b_{i})}.
\end{align*}
\begin{proof}
Clear from the expansion of series and the commutative property of products; similar to adding convergent series.
\end{proof}
\end{prp}

\begin{exa}\label{E: II Root-Two}
\begin{align*}
\prod_{n=1}^{\infty} \bigg[ \frac{(4n-2)^{2}(4n+1)(4n+3)}{(4n-3)(4n-1)(4n+2)^{2}}\bigg]^{n} &= \sqrt{2}
\end{align*}
 \begin{proof}
 This can be proven via the same method outlined in Example \ref{E: II Exponential Catalan}, but instead we will use Proposition \ref{P: Double Product}.
 \begin{align*}
 \prod_{n=1}^{\infty} \bigg[ \frac{(4n-2)^{2}(4n+1)(4n+3)}{(4n-3)(4n-1)(4n+2)^{2}}\bigg]^{n} &=  \prod_{n=1}^{\infty} \bigg[ \frac{(n-1/2)^{2}(n+1/4)(n+3/4)}{(n-3/4)(n-1/4)(n+1/2)^{2}}\bigg]^{n} \\
 &=\prod_{r=0}^{\infty} \prod_{n=1}^{\infty} \frac{(n+r-1/2)^{2}(n+r+1/4)(n+r+3/4)}{(n+r-3/4)(n+r-1/4)(n+r+1/2)^{2}}
 \end{align*}
 By Theorem \ref{T: Gamma Product} we can simplify this back down to a single infinite product,
 \begin{align*}
\prod_{r=0}^{\infty} \prod_{n=1}^{\infty} \frac{(n+r-1/2)^{2}(n+r+1/4)(n+r+3/4)}{(n+r-3/4)(n+r-1/4)(n+r+1/2)^{2}}&=\prod_{r=0}^{\infty}\frac{\Gamma(r+1/4)\Gamma(r+3/4)\Gamma(r+3/2)^{2}}{\Gamma(r+1/2)^{2}\Gamma(r+5/4)\Gamma(r+7/4)}.
\end{align*} 
By the functional equation of the Gamma function, Definition \ref{D: Gamma Function}, and lots of cancelling we arrive at the familiar infinite product found by Catalan\cite{article:7} (the $k=1$ case of Proposition \ref{P: Root Two}).
\begin{align*}
\prod_{r=0}^{\infty}\frac{\Gamma(r+1/4)\Gamma(r+3/4)\Gamma(r+3/2)^{2}}{\Gamma(r+1/2)^{2}\Gamma(r+5/4)\Gamma(r+7/4)} = \prod_{r=0}^{\infty} \frac{(r+1/2)^{2}}{(r+1/4)(r+3/4)} = \prod_{r=0}^{\infty} \frac{(4r+2)^{2}}{(4r+1)(4r+3)}
\end{align*}
\end{proof}
\end{exa}

\begin{prp}\label{P: Gamma Barnes Analogue Finder}
Let $\lbrace a_{1}, ..., a_{i}, b_{1}, ..., b_{i} \rbrace \in\mathbb{C}$ and $a_{1}+...+a_{i}=b_{1}+...+b_{i}$. Then,
\begin{align*}
\prod_{n=1}^{\infty} \frac{(n+a_{1})...(n+a_{i})}{(n+b_{1})...(n+b_{i})}  =  \prod^{\infty}_{n=1}\bigg[ \frac{(n+a_{1})...(n+a_{i})(n+b_{1}+1)...(n+b_{i}+1)}{(n+a_{1}+1)...(n+a_{i}+1)(n+b_{1})...(n+b_{i})}\bigg]^{n}.
\end{align*}
\begin{proof}
Notice, due to the functional equation of the Barnes G-function, Definition \ref{D: Barnes G-Function}, the following formula holds for all complex numbers except for non-positive integers. 
\begin{align*}
\Gamma(z) = \frac{G(z+1)}{G(z)}
\end{align*}
Hence we have a different way to evaluate a product such as those derived from Theorem \ref{T: Gamma Product}. 
\begin{align*}
 \prod_{n=1}^{\infty} \frac{(n+a_{1})...(n+a_{i})}{(n+b_{1})...(n+b_{i})}  =  \frac{G(a_{1})...G(a_{i})}{G(a_{1}+1)...G(a_{i}+1)} \frac{G(b_{1}+1)...G(b_{i}+1)}{G(b_{1})...G(b_{i})}
\end{align*}
All that's left to show is that the right hand side can be written as an the infinite product. Well, this is simply a Type-II product with the constants above. So all we must show is the constants meeting the criteria for Theorem \ref{T: II Product}. \begin{align*}
\sum_{k=1}^{i}a_{i}^{2} + \sum_{k=1}^{i}(b_{i}+1)^{2} &= \sum_{k=1}^{i}a_{i}^{2} + \sum_{k=1}^{i}b_{i}^{2}+2\sum_{k=1}^{i}b_{i}+\sum_{k=1}^{i}1\\
&= \sum_{k=1}^{i}b_{i}^{2} + \sum_{k=1}^{i}(a_{i}+1)^{2}
\end{align*}
The constants meet the requirements hence the proof is complete.
 \end{proof}\end{prp}

\begin{exa}\label{E: II Wallis Product Analogue}
\begin{align*}
\prod_{n=1}^{\infty}\left[ \frac{(2n)^{2}(2n+3)}{(2n-1)(2n+2)^{2}}\right]^{n} = \frac{\pi}{2}.
\end{align*} \begin{proof}
By Proposition \ref{P: Gamma Barnes Analogue Finder}, with $i=2$ and constants
\begin{align*}
a_{1}=a_{2}=0,\quad b_{1}=-\frac{1}{2},\quad b_{2}=\frac{1}{2}.
\end{align*}
We find the equality,
\begin{align*}
\prod_{n=1}^{\infty}\frac{n^{2}}{(n-1/2)(n+1/2)} &= \prod_{n=1}^{\infty}\left[ \frac{n^{2}(n+1/2)(n+3/2)}{(n+1)^{2}(n-1/2)(n+1/2)}\right]^{n}.
\end{align*}
The left-hand side of which is the famous Wallis Product from Example \ref{E: I Wallis Product}, hence we have the Type-II analogue. 
\end{proof}

\end{exa}

\begin{obs} \label{O: Type II Irreducible}
Proposition \ref{P: Gamma Barnes Analogue Finder} implies that for each Type I product we have a Type-II product analogue. One could also ask if the same is true in reverse. If so, we could prove all Type-II products by relating them to their Type-I analogue such as in Example \ref{E: II Root-Two}. If we try the same with Example \ref{E: II Exponential Catalan} the following occurs, \begin{align*}
\prod_{n=1}^{\infty}\bigg[\frac{(4n-1)^{3}(4n+3)}{(4n-3)(4n+1)^{3}}\bigg]^{n} &= \prod_{n=1}^{\infty}\bigg[\frac{(n-1/4)^{3}(n+3/4)}{(n-3/4)(n+1/4)^{3}}\bigg]^{n} \\
&=\prod_{r=0}^{\infty} \prod_{n=1}^{\infty} \frac{(n+r-1/4)^{3}(n+r+3/4)}{(n+r-3/4)(n+r+1/4)^{3}}\\
&=\prod_{r=0}^{\infty} \frac{\Gamma(r+1/4)\Gamma(r+5/4)^{3}}{\Gamma(r+3/4)^{3}\Gamma(r+7/4)}
\end{align*} The expression in the last infinite product can only be written in terms of the Gamma Function (or more complex functions), this implies there are some Type-II products that cannot be written as Type-I products. 
\end{obs}

\begin{obs}\label{O: Irrationality}
It should be noted that Type-II products can converge to many numbers of different irrationality. Example \ref{E: II Exponential Catalan} has unknown irrationality, Example \ref{E: II Rational} is rational, Example \ref{E: II Root-Two} is algebraic irrational and Example \ref{E: II Wallis Product Analogue} is transcendental.
\end{obs}

\section{Multiple Gamma Function}

\begin{thm} \label{T: N product}
Let $n\in\mathbb{N}_{+}$ and suppose that there exists finite sequences $(a_{k})_{k=1}^{i}$ and $(b_{k})_{k=1}^{i}$ with $a_{k}$, $b_{k}\in\mathbb{C}\setminus\mathbb{Z}_{<0}$, such that $\sum_{k=1}^{i} a_{k}^{j} = \sum_{k=1}^{i} b_{k}^{j}$, for $j\in\lbrace 1,2,...,n\rbrace$. Then, \begin{align*}
\prod_{k=1}^{\infty}\left[ \frac{(k+a_{1})...(k+a_{i})}{(k+b_{1})...(k+b_{i})}\right]^{{n+k-2}\choose{n-1}} = \frac{\Gamma_{n}(b_{1}+1)...\Gamma_{n}(b_{i}+1)}{\Gamma_{n}(a_{1}+1)...\Gamma_{n}(a_{i}+1)}.
\end{align*} \emph{Products of this form will be denoted as Type-$n$ products.}
\begin{proof}
The Multiple Gamma function can be written explicitly as follows \cite{ARTICLE:2},
\begin{align*}
\Gamma_{n}(z+1) &= e^{Q_{n}(z)}\prod_{k=1}^{\infty}\left\lbrace \left( 1+\frac{z}{k} \right)^{-{{n+k-2}\choose{n-1}}}\exp\left[ {{n+k-2}\choose{n-1}}\left( \sum_{j=1}^{n}\frac{(-1)^{j-1}}{j}\frac{z^{j}}{k^{j}} \right) \right] \right\rbrace, \\
Q_{n}(z) &= (-1)^{n-1}\left[ -zA_{n}(1) + \sum_{k=1}^{n-1}\frac{p_{k}(z)}{k!}\left(f_{n-1}^{(k)}(0) - A_{n}^{(k)}(1) \right) \right], \\
f_{n}(z) &=  -zA_{n}(1) + \sum_{k=1}^{n-1}\frac{p_{k}(z)}{k!}\left(f_{n-1}^{(k)}(0) - A_{n}^{(k)}(1) \right) + A_{n}(z),  \\
A_{n}(z) &= \sum_{k=1}^{\infty} (-1)^{n-1}{{n+k-2}\choose{n-1}}\left[ -\ln\left( 1 + \frac{z}{k}\right) +   \sum_{j=1}^{n}\frac{(-1)^{j-1}}{j}\frac{z^{j}}{k^{j}} \right],\\
p_{n}(z) &= \frac{1}{n+1} \sum_{k=1}^{n+1} {{n+1}\choose{k}}B_{n+1-k}z^{k},\quad (n\in\mathbb{N}).
\end{align*} In the same way we proved Theorem \ref{T: II Product} we can multiply multiple instances of these infinite products for the Multiple Gamma Function together. So consider $(a_{k})$,$(b_{k}) \in\mathbb{C}\setminus\mathbb{Z}_{<0}$. Then,
\begin{align*}
\frac{\Gamma_{n}(b_{1}+1)...\Gamma_{n}(b_{i}+1)}{\Gamma_{n}(a_{1}+1)...\Gamma_{n}(a_{i}+1)} = e^{S}\prod_{k=1}^{\infty}\left\lbrace \left[ \frac{(k+a_{1})...(k+a_{i})}{(k+b_{1})...(k+b_{i})}\right]^{{n+k-2}\choose{n-1}} \exp\left[ {{n+k-2}\choose{n-1}}\left( \sum_{j=1}^{n}T_{j} \right) \right] \right\rbrace 
\end{align*}\vspace{-0.3cm} \begin{align*}
\text{where}\quad T_{j} &= \frac{(-1)^{j-1}}{jk^{j}}\Big( \big(b_{1}^{j}+b_{2}^{j}+ ... + b_{i}^{j}\big) - \big(a_{1}^{j}+a_{2}^{j}+ ... + a_{i}^{j}\big)\Big), \\
S &= \big(Q_{n}(b_{1}) + Q_{n}(b_{2}) + ... + Q_{n}(b_{i})\big) - \big( Q_{n}(a_{1}) + Q_{n}(a_{2}) + ... + Q_{n}(a_{i})\big),
\end{align*} 
(it can be written as one product as each individual product is itself convergent). Now, suppose finite sequences $(a_{k})$, $(b_{k})$, satisfy the theorem criteria. Finite sums commute, so for each $j\in\lbrace 1,2,...,n \rbrace$, $T_{j}=0$. Notice, similar to why $T_{j} = 0$,  \begin{align}
\big(p_{n}(b_{1})+...+p_{n}(b_{i})\big) - \big(p_{n}(a_{1}) +... +p_{n}(a_{i})\big) = 0.
\end{align} And again using a similar argument we see using (1) and the theorem criteria for $j=1$, \begin{align*}
S =\big(Q_{n}(b_{1})+...+Q_{n}(b_{i})\big) - \big(Q_{n}(a_{1}) +... +Q_{n}(a_{i})\big) = 0.
\end{align*} This shows both the exponentials inside and outside the product are $1$, completing the proof.
\end{proof}
\end{thm}

\begin{obs}
Almost all Propositions and Theorems discussed previously can be generalised. They were explored in greater detail before as constants which satisfy the theorem criteria become less frequent for greater $n$. It is easy to see Proposition \ref{P: Double Product} and Proposition \ref{P: Gamma Barnes Analogue Finder} can be extended. Following from this, if Proposition \ref{P: Gamma Barnes Analogue Finder} is extended to any Type-$n$ product then we can elaborate on Remark \ref{O: Irrationality} to say that all Type-$n$ products can converge to constants of different rationality.
\end{obs}

\begin{obs} \label{O: Dup}
For higher values of $n$ it becomes a much more arduous task finding solutions which fit the criteria required. Even when solutions are found, simple closed forms are even harder to find as there is limited research on the Multiple Gamma Function for $n \geq 3$. For example, the Barnes G-Function has a duplication formula seen below \cite{ARTICLE:1}, but for higher values of $n$ the formula gets very ugly and the literature is very sparse. \begin{align*}
G(2z)= e^{-\frac{1}{4}}A^{3}2^{2z^{2}-3z+\frac{11}{12}}\pi^{\frac{1}{2}-z}G(z)G\left(z +\frac{1}{2} \right)^{2}G(z+1),
\end{align*} where $A$ is the Glaisher-Kinkelin constant from Definition \ref{D: GK Constant}. Finding solutions to the theorem criteria is equivalent to finding solutions (not necessarily ideal solutions) to the Prouhet-Tarry-Escott problem hence one can use existing research on this open problem.
\end{obs}

\setcounter{equation}{0}

\begin{exa}
\begin{align*}
\prod_{n=1}^{\infty}\left[ \frac{(2n)^{3}(2n+3)^{3}}{(2n-1)(2n+1)^{2}(2n+2)^{2}(2n+4)} \right]^{n^{2}} = \frac{8A^{12}}{e\pi^{3}\sqrt[3]{2}}
\end{align*}
\begin{proof}
Using Theorem \ref{T: N product} for $n=3$, $i=6$ and with the following constants, \begin{align*}
a_{1}=a_{2}&=a_{3}=0, \quad a_{4}=a_{5}=a_{6}=\frac{3}{2}, \\
b_{1}=-\frac{1}{2},\quad  b_{2}&=b_{3}=\frac{1}{2}, \quad b_{4}=b_{5}=1, \quad b_{6} = 2.
\end{align*} One can easily check these satisfy the theorem criteria. Then we have,
\begin{align}
\prod_{n=1}^{\infty}\left[ \frac{(2n)^{3}(2n+3)^{3}}{(2n-1)(2n+1)^{2}(2n+2)^{2}(2n+4)} \right]^{\frac{1}{2}n(n+1)} = \frac{\Gamma_{3}(1/2)\Gamma_{3}(3/2)^{2}\Gamma_{3}(2)^{2}\Gamma_{3}(3)}{\Gamma_{3}(1)^{3}\Gamma_{3}(5/2)^{3}}.
\end{align} Then we can use the functional equation of the Multiple Gamma Function from Definition \ref{D: Multiple Gamma Function} to simplify, \begin{align*}
\frac{\Gamma_{3}(1/2)\Gamma_{3}(3/2)^{2}\Gamma_{3}(2)^{2}\Gamma_{3}(3)}{\Gamma_{3}(1)^{3}\Gamma_{3}(5/2)^{3}} = \frac{1}{\pi^{3/2}G(1/2)^{4}}.
\end{align*} Then, using Theorem \ref{T: II Product} with the same constants (using a similar manipulation method) we also have, \begin{align}
\prod_{n=1}^{\infty}\left[ \frac{(2n)^{3}(2n+3)^{3}}{(2n-1)(2n+1)^{2}(2n+2)^{2}(2n+4)} \right]^{n} = \frac{\pi^{2}}{8}.
\end{align} A formula is known for $G(1/2)$ in terms of famous mathematical constants \cite{article:6}, (this also follows from the duplication formula from Remark \ref{O: Dup})
\begin{align*}
G(1/2) = \frac{2^{1/24}e^{1/8}}{A^{3/2}\pi^{1/4}}.
\end{align*} Hence squaring (1) and dividing by (2) gives the required result.
\end{proof}
\end{exa}

\section*{Conclusion}
We have seen a more systematic approach to finding infinite products similar to that of Wallis' Product, this lead to results involving nested roots and $\varphi$. We have generalised upon that to the Barnes G-function, where we have proved several propositions relating to finding analogues of Wallis-type products; from which we produced some beautiful closed form representations of Type-II products. This again can be further generalised to the Multiple Gamma Functions, but due to restrictions in theorem criteria, Type-$n$ products (for $n\geq 3$) are much more difficult to find. We gave one example independent from previous theorems and also remarked on the potential of the simpler propositions. 
 \pagebreak

\bibliography{bib4}

\end{document}